\documentclass[a4paper,11pt]{article}
\usepackage[margin=1cm, bottom = 2cm]{geometry}
\usepackage[russian]{babel}
\usepackage{multicol}
\usepackage{amsmath}
\usepackage{graphicx}
\usepackage{caption}
\renewcommand{\section}[2]{}
\title{Синхронизация математических моделей \\нейронной активности Моррис-Лекара\footnote{Эта работа поддержана грантом СПбГУ ID 84912397.}}
\author{А. В. Рыбалко, Д. М. Семёнов, А. Л. Фрадков\\
Кафедра Теоретической кибернетики СПбГУ, Санкт-Петербург, Россия\\
Институт проблем машиноведения РАН, Санкт-Петербург, Россия}
\date{}
\begin{document}
\maketitle

Данная работа посвящена задаче синхронизации двух моделей нейронов Моррис-Лекара. Модель Моррис-Лекара представляет из себя систему дифференциальных уравнений второго порядка, описывающую непростую взаимосвязь между мембранным потенциалом и активацией ионных каналов внутри мембраны. Синхронизация, то есть такое состояние сети, при котором модели начинают вести себя одинково в некотором смысле, такого рода моделей представляет не только математический интерес, но и биологический: синхронная активность нейронных популяций играет очень важную роль в работе мозга. Для решения этой задачи был применён алгоритм скоростного градиента, являющегося непрерывным аналогом градиентных алгоритмов. С помощью него был получен алгоритм управления силой связи нейронов, обеспечивающий достижение поставленной цели. Моделирование в среде MATLAB продемонстрировало корректность полученного результата и быструю сходимость соответствующих переменных моделей нейронов друг к другу.

\begin{multicols}{2}[\columnsep=0.8cm]
\begin{center}
    \textbf{Введение}
\end{center}

Большой интерес представляет изучение синхронной активности нейронов. Установлено, что во-первых, синхронизированная активность больших популяций нейронов является главным механизмом образования ритмов головного мозга человека, играющих важную роль в процессах передачи и обработки информации в центральной нервной системе; а во-вторых, чрезмерная локальная синхронизация в некоторых сетях имеет патологические проявления, такие как тремор в болезни Паркинсона, эпилептические припадки, а также некоторые расстройства высшей мозговой функции (шизофрения, аутизм, и т. д.). Разработка методов для контроля за невозникновением нежелательной синхронной мозговой активности является важнейшей клинической задачей.

Сущестование математических моделей, описывающих активность нейронов с помощью систем дифференциальных уравнений, даёт возможность математически исследовать явление синхронизации, что позволит улучшать результаты лечения и уменьшать его побочные эффекты. Биологический нейрон можно моделировать на различных уровнях абстракции, но суть многих моделей заключается в изображении нервной клетки как осциллятора определённого уровня сложности. Эта работа посвящена поиску условия, при котором два идентичных нейрона полностью синхронизируются. Для этого была выбрана модель Моррис-Лекара \cite{morrislecar}, состоящая из двух дифференциальных уравнений и являющаяся упрощённой моделью Ходжкина-Хаксли. Она включает в себя кинетику калиевых и кальциевых ионных каналов, представляющих собой сложные белковые структуры с молекулярными системами открытия, закрытия,
селективности, инактивации и регуляции. Каналы так же принимают участие в процессах передачи информации с одной нервной клетки на другую, обеспечивают создание мембранного потенциала покоя, возбудимость, а также активную или пассивную деполяризацию (понижениие потенциала покоя), инициируют выделение гормонов и сокращение мышечных волокон.

\begin{center}
    \textbf{Построение алгоритма синхронизации}
\end{center}

Математическая модель нейронной активности Моррис-Лекара имеет вид \cite{morrislecar}:
\begin{equation}
\begin{cases}
\begin{split}
   \dot V =\, & \frac{1}{C}(- g_L(V - V_L)
   - g_{Ca}m(V - V_{Ca})  \\&\quad\quad\quad\quad- g_KN(V - V_K) + I),\\
   \dot N =\, & \lambda\frac{n - N}{\tau},
\end{split}
\end{cases} \label{ML}
\end{equation}

$m = \frac{1}{2}(1 + \tanh(\frac{V - \widetilde{V_1}}{\widetilde{V_2}})),$

$n = \frac{1}{2}(1 + \tanh(\frac{V - \widetilde{V_3}}{\widetilde{V_4}})).$

$\tau = \frac{1}{\cosh(\frac{V - \widetilde{V_3}}{2\widetilde{V_4}})},$\\
где

\begin{itemize}
  \item V - мембранный потенциал
  \item N - доля открытых каналов \begin{math}K^{+}\end{math} 
  \item С - мембранная ёмкость 
  \item \begin{math}g_{Ca}, g_K, g_L\end{math} - максимальные проводимости трансмембранных токов для \begin{math}Ca^{+}, K^{+}\end{math} и утечки 
  \item \begin{math}V_{Ca}, V_K, V_L\end{math} - реверсивные потенциалы, соответствующие \begin{math}Ca^{+}, K^{+}\end{math} и утечке 
  \item I - внешний ток
  \item \begin{math}\widetilde{V_1}, \widetilde{V_2}, \widetilde{V_3}, \widetilde{V_4}\end{math} - параметры настройки для стационарного состояния
  \item \begin{math}\lambda\end{math} - масштабирующий парамер
\end{itemize}

Теперь перейдём в поиску условия синхронизации двух нейронов. Для этого сначала определим, что в данной работе понимается под синхронизацией: синхронизацией будем называть согласованное во времени функционирование двух или нескольких процессов или объектов \cite{cyb_ph}. Или на математическом языке:

\textbf{Определение.} Будем говорить, что имеет место синхронизация процессов \begin{math}x^{(i)}(t), i = 1, 2,...,k,\end{math} относительно характеристики \begin{math}C_t\end{math} и функции сравнения \begin{math}F_i\end{math}, если существуют вещественные числа (временные или фазовые сдвиги) \begin{math}\tau_i, i = 1,2,...,k\end{math} такие, что для всех \begin{math}t\geq0\end{math} выполняются соотношения $F_1(C_{t+\tau_1}[x_1]) = F_2(C_{t+\tau_2}[x_2]) = ... = F_k(C_{t+\tau_k}[x_k]).$

Мы будем использовать вид синхронизации, называемый координатным. Координатная синхронизация укладывается в предложенное выше общее определение, если положить \begin{math}C_t(x_i) = x_i(t)\end{math}, где через \begin{math}x_i(t)\end{math} обозначено значение вектора состояния i-й подсистемы в момент времени t, а функции сравнения взять тождественными: \begin{math}F_i(x)=x, i = 1,..., k.\end{math}

Рассмотрим две связанные модели Моррис-Лекара и найдём условие, при котором они будут синхронизироваться. В рамках этой работы будем исследовать модели с одинковыми наборами параметров, различающиеся только начальными данными:

\begin{equation}
\begin{cases}
\begin{split}
   \dot V_1 =\, & \frac{1}{C}(- g_L(V_1 - V_L)
   - g_{Ca}m_1(V_1 - V_{Ca})  \\&- g_KN_1(V_1 - V_K) + I) + \sigma(V_2 - V_1),\\
   \dot N_1 =\, & \lambda\frac{n_1 - N_1}{\tau_1},\\
   \dot V_2 =\, & \frac{1}{C}(- g_L(V_2 - V_L)
   - g_{Ca}m_2(V_2 - V_{Ca})  \\&- g_KN_2(V_2 - V_K) + I) + \sigma(V_1 - V_2),\\
   \dot N_2 =\, & \lambda\frac{n_2 - N_2}{\tau_2},
\end{split}
\end{cases}
\end{equation}

$m_i = \frac{1}{2}(1 + \tanh(\frac{V_i - \widetilde{V_1}}{\widetilde{V_2}})),$

$n_i = \frac{1}{2}(1 + \tanh(\frac{V_i - \widetilde{V_3}}{\widetilde{V_4}})).$

$\tau_i = \frac{1}{\cosh(\frac{V_i - \widetilde{V_3}}{2\widetilde{V_4}})},$\\
где $i = 1,\,2$, а \begin{math}\sigma(t)\end{math} обозначает силу связи между нейронами.

Или перепишем в упрощённой форме:

\begin{equation}
\begin{cases}
\begin{split}
   \dot V_1 &=\, f(V_1, N_1) + \sigma(V_2 - V_1),\\
   \dot N_1 &=\, g(V_1, N_1),\\
   \dot V_2 &=\, f(V_2, N_2) + \sigma(V_1 - V_2),\\
   \dot N_2 &=\, g(V_2, N_2),
\end{split}
\end{cases} \label{ML_simpl}
\end{equation}

Теперь перейдём к новым переменным \begin{math}e_V := V_1 - V_2, e_N := N_1 - N_2\end{math}:
\begin{equation}
 \begin{cases}
   \dot e_V = F(e_V, e_N, \psi_1(t)) - 2\sigma e_V,\\
   \dot e_N = G(e_V, e_N, \psi_2(t)),
 \end{cases} \label{ML_diff}
\end{equation}
где $F(e_V, e_N, \psi_1(t)) = f(V_1, N_1) - f(V_2, N_2),$ $ G(e_V, e_N, \psi_2(t)) = g(V_1, N_1) - g(V_2, N_2) $, a $\psi_1(t), \psi_2(t)$ "--- нелинейности, явно не выражающиеся через переменные $e_V$ и $e_N$. 

Итак, цель управления для достижения синхронизации, исходя из определения данного выше, "--- это найти условие, при котором
\begin{equation}
e_V \rightarrow 0,\,
e_N \rightarrow 0 \,\,\mbox{при}\,\, t \rightarrow +\infty
\end{equation}
Иными словами, чтобы гарантировать синхронизацию в системе \eqref{ML_simpl}, необходимо обеспечить устойчивость системы \eqref{ML_diff} \cite{semenov}. Мы будем делать это через управление, определяемое функцией $\sigma(t)$. Чтобы не возникало неопределённостей, будем считать, что \begin{math}|V_1| \leq M_1, |V_2| \leq M_2, |N_1| \leq M_3, |N_2| \leq M_4 \ \forall t \geq t_0 + T\end{math} для некоторого T и некоторых \begin{math}M_1,...,M_4\end{math}. Такие оценки естественны и соответствуют биологическому смыслу переменных мембранного потенциала и доли открытых каналов $K^{+}$. 

Далее будем использовать метод скоростного градиента. Этот метод предназначен для решения задач управления непрерывными по времени системами, в которых цель управления задана при помощи гладкой целевой функции \cite{fradkov90}. В качестве такой целевой функции возьмем $Q(e(t)) = \frac{1}{2}(e_V^2+e_N^2)$, тогда цель усправления можно сформулировать так: $Q(e(t)) \rightarrow 0$ при $t \rightarrow \infty$, основываясь на следующих соображениях. Как уже было замечено ранее, нам необходимо обеспечить устойчивость системы \eqref{ML_diff}, для этого можно использовать 1 теорему Ляпунова \cite{miroshnik}: целевая функция $Q(e(t))$ подходит под определение функции Ляпунова, а производная данной функции Ляпунова в силу системы \eqref{ML_diff} будет неположительной, начиная с некторого момента времени, как раз когда $Q(e(t)) \rightarrow 0$ при $t \rightarrow \infty$, то есть убывает, учитывая положительную определённость фунции $Q(e(t))$.

В методе скоростного градиента сначала вычисляется производная в силу системы \eqref{ML_diff}: 
\begin{equation}
\begin{split}
\omega(e, \sigma, t)& := \dot Q(e(t)) = e_V\dot e_V + e_N\dot e_N \\= e_V&(F(e_V, e_N, \psi_1(t)) - 2\sigma e_V) \\&+ e_N G(e_V, e_N, \psi_2(t)).
\end{split}
\end{equation}

Затем находится градиент полученной функции $\dot Q(e(t))$ по входной переменной $\sigma$: 
\begin{equation}
\nabla_\sigma \ \omega(e, \sigma, t) = -2e_V^2.
\end{equation}

Наконец, дифференциальным уравнением задаётся алгоритм изменения $\sigma(t)$: 
\begin{equation}
    \dot \sigma = - \gamma \nabla_\sigma \ \omega(e, \sigma, t) = 2\gamma e_V^2,
\end{equation}
где $\gamma$ > 0 - коэффициент усиления. 

В результате, получаем новую систему дифференциальных уравнений, моделирующую динамику двух синхронизированных нервных клеток: 
\begin{equation}
 \begin{cases}
   \dot V_1 = f(V_1, N_1) + \sigma(V_2 - V_1)\\
   \dot N_1 = g(V_1, N_1)\\
   \dot V_2 = f(V_2, N_2) + \sigma(V_1 - V_2)\\
   \dot N_2 = g(V_2, N_2)\\
   \dot \sigma = 2\gamma (V_1 - V_2)^2.
 \end{cases} \label{fin_syst}
\end{equation} 

\begin{center}
    \textbf{Результаты компьютерного моделирования}
\end{center}

На Рис. 1 и Рис. 2 представлено поведение переменных мембранного потенциала и доли открытых каналов \begin{math}K^{+}\end{math} соответственно для некоторого набора параметров и начальных данных \begin{math}V_0, N_0\end{math}: 
\begin{equation}
\begin{split}
&V_0 = -35, N_0 = 0.9, \\&\widetilde{V_1} = -1, \widetilde{V_2} = 15, \widetilde{V_3} = 10, \widetilde{V_4} = 14.5, \\&g_{Ca} = 4, g_K = 8, g_L = 2, \\&V_{Ca} = 100, V_K = -70, V_L = -50, \\&C = 20, I = 50, \lambda = 0.1. 
\end{split} \label{init1}
\end{equation}

На Рис. 3 и Рис. 4 представлено поведение переменных мембранного потенциала и доли открытых каналов $K^+$ системы  для некоторого набора начальных данных и параметров системы \eqref{fin_syst}:
\begin{equation}
\begin{split}
    &V_{10} = -35, N_{10} = 0.9, \\&V_{20} = 10, N_{20} = 3, \sigma_0 = -1, \\&\widetilde{V_1} = -1, \widetilde{V_2} = 15, \widetilde{V_3} = 10, \widetilde{V_4} = 14.5, \\&g_{Ca} = 4, g_K = 8, g_L = 2, \\&V_{Ca} = 100, V_K = -70, V_L = -50, \\&C = 20, I = 50, \lambda = 0.1, \gamma = 0.5.
\end{split} \label{init2}
\end{equation}

\begin{minipage}{\linewidth}
      \includegraphics[width=0.9\linewidth]{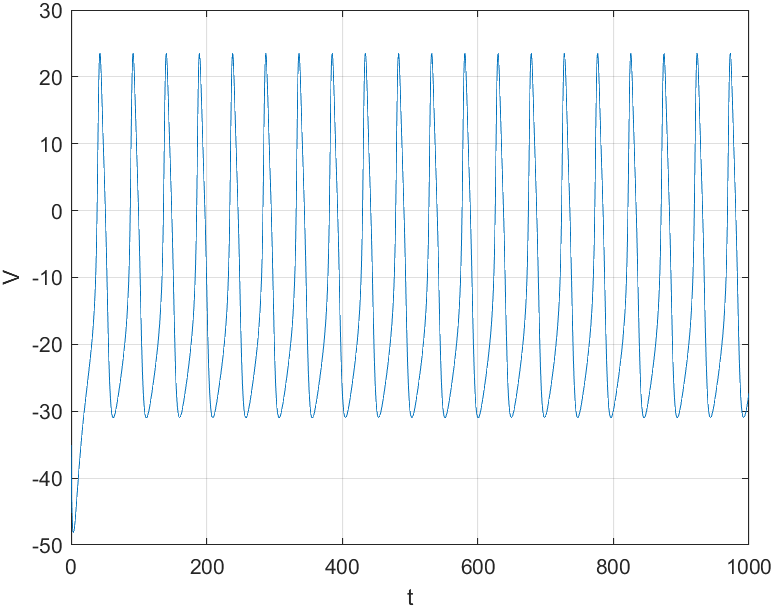}
      \captionof{figure}{График первой переменной решения\\системы \eqref{ML} для набора начальных данных и \\параметров \eqref{init1}.}
\end{minipage}

\begin{minipage}{\linewidth}
      \includegraphics[width=0.9\linewidth]{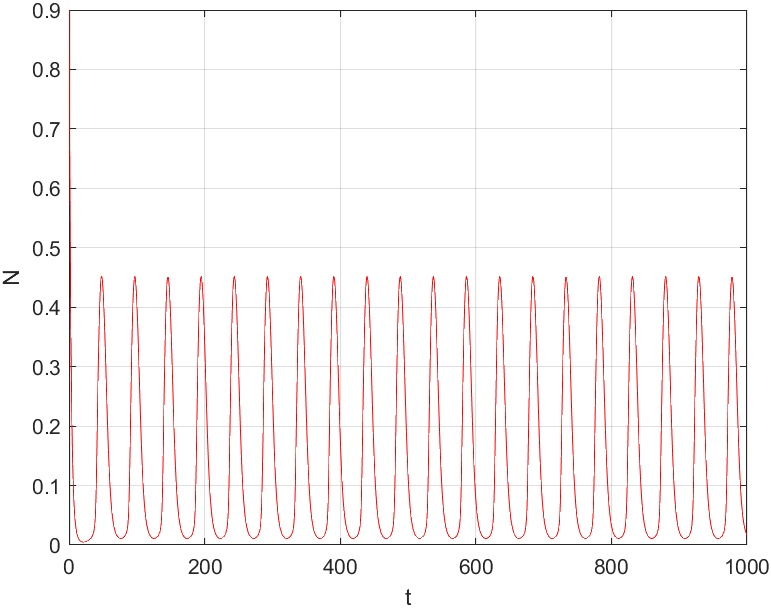}
      \captionof{figure}{График второй переменной решения\\системы \eqref{ML} для набора начальных данных и \\параметров \eqref{init1}.}
\end{minipage} 

\begin{minipage}{\linewidth}
      \includegraphics[width=0.9\linewidth]{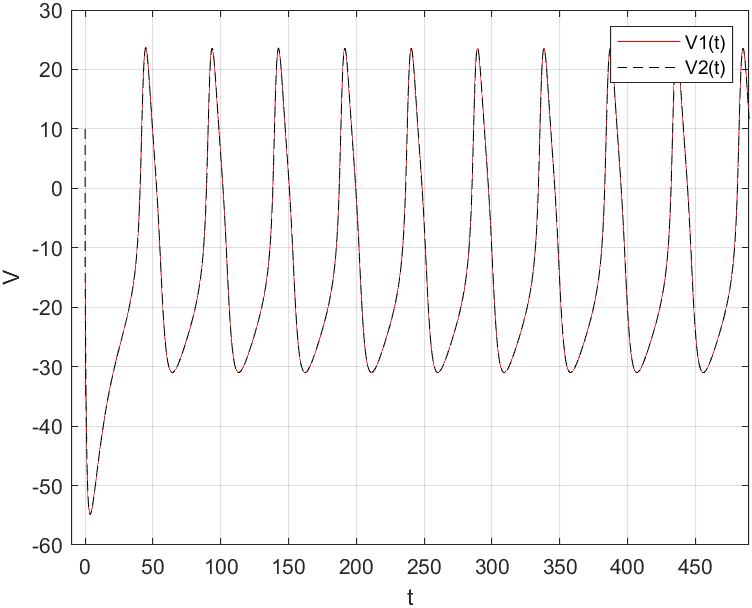}
      \captionof{figure}{График первой переменной решения\\системы \eqref{fin_syst} для набора начальных данных и \\параметров \eqref{init2}.}
\end{minipage}

\begin{minipage}{\linewidth}
      \includegraphics[width=0.9\linewidth]{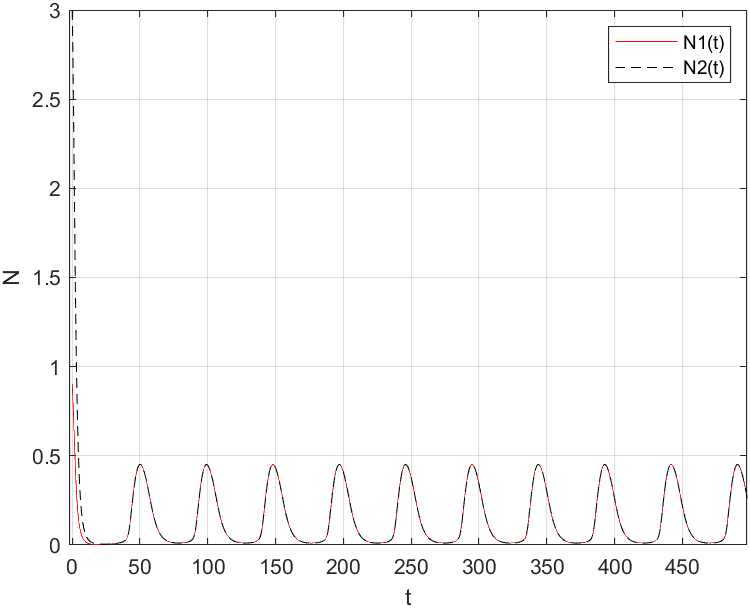}
      \captionof{figure}{График второй переменной решения\\системы \eqref{fin_syst} для набора начальных данных и \\параметров \eqref{init2}.}
\end{minipage}
\vspace{1cm}

\begin{minipage}{\linewidth}
      \includegraphics[width=0.9\linewidth]{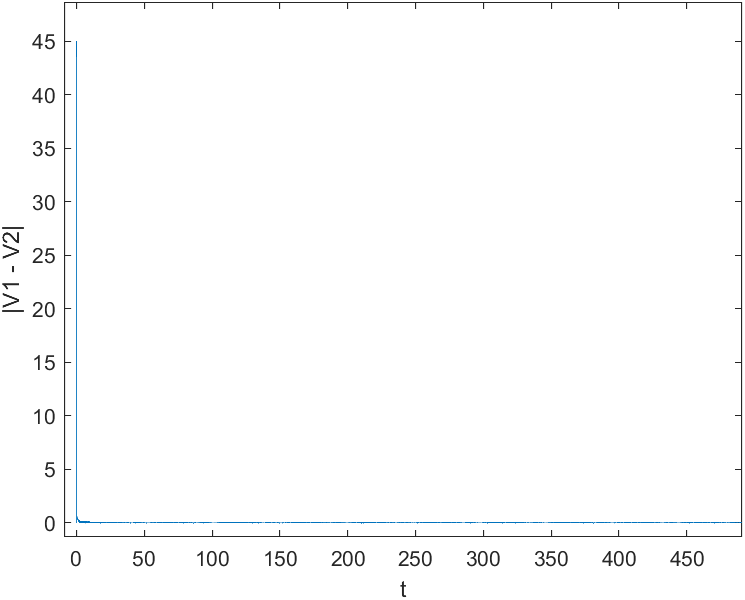}
      \captionof{figure}{График абсолютной разности переменных \\$V_1$ и $V_2$ решения системы \eqref{fin_syst} для набора \\начальных данных и параметров \eqref{init2}.}
\end{minipage}
\vspace{1cm}

\begin{minipage}{\linewidth}
      \includegraphics[width=0.9\linewidth]{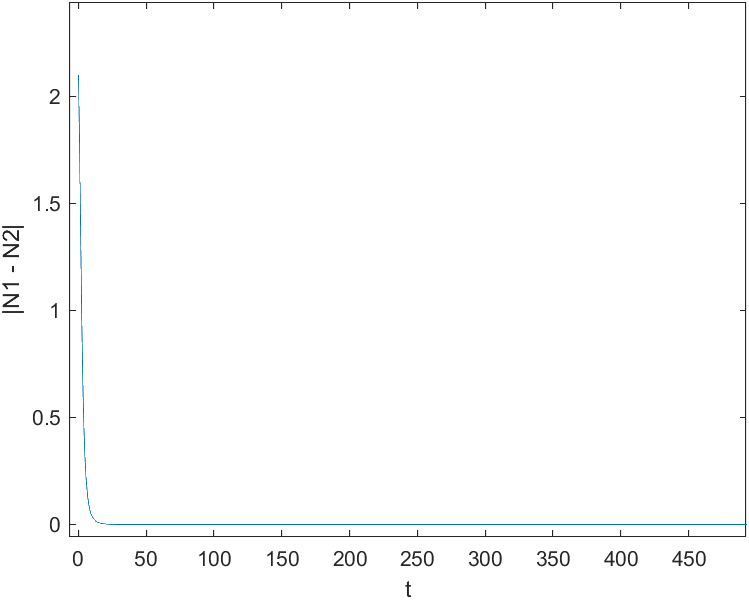}
      \captionof{figure}{График абсолютной разности переменных \\$N_1$ и $N_2$ решения системы \eqref{fin_syst} для набора \\начальных данных и параметров \eqref{init2}.}
\end{minipage}

Для наглядности на Рис. 5 и Рис. 6 изображены абсолютные разности переменных $V_1$ и $V_2$, $N_1$ и $N_2$.

Очевидно, что даже при большой разнице между начальными данными, достигается цель управления за конечное время, иными словами, поведение нейронов очень быстро синхронизируется в опредлённом выше смысле "--- поставленная цель достигнута. 

\begin{center}
    \textbf{Заключение}
\end{center}

В данной работе широко распространённая в нейронауках задача изучения условий синхронизации моделей нейронов решена через управление силой связи нейронов. Адаптивный алгоритм управления для полной синхронизации двух моделей нервных клеток Моррис-Лекара получен с помощью метода скоростного градиента. Полученное решение было промоделировано в среде MATLAB. Результаты моделирования продемонстрировали быстрое и точное достижение цели управления.

\begin{center}
    \textbf{Список литературы}
\end{center}

\vspace{5cm}

\end{multicols}

\newpage

\newpage
\begin{center}
    \LARGE
    Synchronization of Morris-Lecar\\ mathematical models of neural activity\\ 
    \vspace{0.5cm}
    \large
    A. V. Rybalko, D. M. Semenov, A. L. Fradkov\\
Theoretical Cybernetics Department, SPBU, Saint Petersburg, Russia\\
Institute for Problems in Mechanical Engineering, RAS, Saint Petersburg, Russia
\end{center}
\vspace{0.5cm}

\textbf{Acknowledgements:} This paper was supported by grant SPbU ID 84912397.\\

\textbf{Abstract}\\

This work is devoted to the problem of synchronization of two Morris-Lecar neuron models. The Morris-Lecar model is a second-order system of differential equations, which describes an uneasy relationship between the membrane potential and the activation of ion channels inside the membrane. Synchronization, which is a state of a network when models begin to act similarly in some sense, of such models is interesting not only from a mathematical point of view, but also from a biological one, because synchronous activity plays a very important role in brain functioning. The speed gradient algorithm, which is a continuous version of gradient algorithms, was applied to solve this problem. The algorithm of coupling strength control was obtained. It ensures the achievement of the control goal. The MATLAB modelling demonstrated the correctness of the obtained result and fast convergence rate of corresponding models' variables to each other.

\end{document}